\documentstyle[amssymb,amsfonts]{amsart}

\def\Q{{\hbox{\bf Q}}}

\def\be{{\bf{e}}}
 at 10 true pt

\def\be#1{ \begin{equation}\label{#1} }
\def\bas{\begin{align*}}
\def\eas{\end{align*}}
\def\bi{\begin{itemize}}
\def\ei{\end{itemize}}

\def\Z{{\hbox{\bf Z}}}

\newenvironment{proof}{\noindent {\bf Proof} }{\endprf\par}
\def \endprf{\hfill  {\vrule height6pt width6pt depth0pt}\medskip}
\def\emph#1{{\it #1}}
\def\textbf#1{{\bf #1}}

\def\BZ{{\mathbf Z}}

\def\mod {\hbox {\rm mod}}

\theoremstyle{plain}
  \newtheorem{theorem}[subsection]{Theorem}
  \newtheorem{question}[subsection]{Question}

  \newtheorem{proposition}[subsubsection]{Proposition}
  \newtheorem{claim}[subsubsection]{Claim}
  \newtheorem{lemma}[subsection]{Lemma}
  \newtheorem{corollary}[subsection]{Corollary}

\theoremstyle{remark}

\theoremstyle{definition}

\include{psfig}

\begin{document}
\title[Classification theorems for sumsets modulo a prime]
{Classification theorems for sumsets modulo a prime}

\author{Hoi H. Nguyen}
\address{Department of Mathematics, Rutgers University, 110 Frelinghuysen rd., Piscataway, NJ 08854, USA}
\email{hoi@@math.rutgers.edu}
\thanks{}

\author{Van H. Vu}
\address{Department of Mathematics, Rutgers University, 110 Frelinghuysen rd., Piscataway, NJ 08854, USA}
\email{vanvu@@math.rutgers.edu}

\thanks{This work was written while the first author was supported by a DIMACS summer research fellowship, 2006.}

\thanks{The second author is an A. Sloan Fellow and is supported by an NSF Career Grant.}

\begin{abstract}
Let $\BZ_p$ be the finite field of prime order $p$ and $A$ be a
subsequence of $\BZ_p$. We prove several classification results
about the following questions:

(1) When can one represent zero as a sum of some elements of $A$ ?

(2) When can one represent every element of $\BZ_p$ as a sum of some
elements of $A$ ?

(3) When can one represent every element of $\BZ_p$ as a sum of $l$
elements of $A$ ?
\end{abstract}

\maketitle
\section{Introduction.}

\noindent Let $G$ be an additive group and $A$ be a sequence of (not
necessarily different) elements of $G$. We denote by $S_A$ the
collection of partial sums of $A$

$$
\sum(A) := \left \{ \sum_{x\in B}x \ |\ \emptyset \neq B\subset A,
|B| < \infty \right \}.
$$

\noindent For a positive integer $l\le |A|$ we denote by
$\sum_{l}(A)$ the collection of partial sums of $l$ elements of $A$,

$$
\sum_l(A) := \left \{ \sum_{x\in B}x \ | \ B\subset A, |B|=l \right
\}.
$$

{\bf \noindent Example.} If $G= \BZ_{11}, A =\{1,1,7 \}$ then
$\sum(A)= \{1,2,7,8,9 \}$ and $\sum_2(A)=\{2,8\}$.

\vskip2mm

\noindent The following questions are among the most popular  in
Additive Combinatorics.

\begin{question} \label{question:1} When is $ 0 \in \sum(A)$ and when is $\sum(A)= G$ ?
\end{question}

\begin{question} \label{question:2} For a given $l$ when is $0 \in \sum_l(A)$ and when is $\sum_l(A)=G$?
\end{question}

\noindent There is a vast amount of results concerning these
questions (see for instance \cite{GG},\cite{Ge},\cite{Sun}),
including classical results such as Olson's theorem and the
Erd\H{o}s-Ginzburg-Ziv theorem.

\noindent If $ 0 \notin \sum(A)$ (or, respectively, $0 \notin
\sum_{l}(A)$), then we say that $A$ is {\it zero-sum-free} (or,
respectively, $l$-{\it zero-sum-free}). If $\sum(A) = G$ (or,
respectively, $\sum_l(A) = G$), then we say that $A$ is {\it
complete} (or, respectively, $l$-{\it complete}); and otherwise we
say that $A$ is {\it incomplete} ($l$-{\it incomplete}).

\noindent We will focus on the case $G=\BZ_p$, the cyclic group of
order $p$, where $p$ is a large prime. The main goal of this paper
is to give a strong classification for zero-sum-free, incomplete and
$l$-incomplete sequences of $\Z_p$. These classifications refine and
extend an implicit result in \cite{SzemVu1}. Together they support
the following general phenomenon:

\vskip2mm

{\it The main reason for a  sequence to be  zero-sum-free or
incomplete is that its elements have small norm. }

\vskip2mm

\noindent For instance, if the elements of a sequence (viewed as
positive integers between $0$ and $p-1$) add up to a number less
than $p$, then the sequence is clearly zero-sum-free. One of our
results, Theorem \ref{theorem:main1}, shows that any zero-sum-free
sequence in $\Z_p$ can be brought into this form after a dilation
and after truncation of a negligible subset.

\noindent Our results  have many applications (see Sections
\ref{section:zerofree},\ref{section:incomplete},\ref{section:partition}
and \ref{section:l-incomplete}). In particular,  we will prove a
refinement of the well-known Erd\H{o}s-Ginzburg-Ziv theorem (see
Section \ref{section:l-incomplete}). The common theme of  these
applications is the following.

\vskip2mm

{\it Any long zero-sum-free or incomplete sequence is a subsequence
of a unique extremal sequence (after a proper linear transformation
and a possible truncation of a negligible subsequence). }

\vskip2mm

\noindent In the rest of this section, we introduce our notation.
The remaining sections are organized as follows. In Section
\ref{section:result}, we present our classification theorems.
Sections
\ref{section:zerofree},\ref{section:incomplete},\ref{section:partition},\ref{section:l-incomplete}
are devoted to applications. Section \ref{section:lemma} contains
the main lemmas needed for the proofs. The proofs of the
classification theorems come in Sections
\ref{section:main3},\ref{section:main2} and \ref{section:main1}.

\vskip2mm

{\bf \noindent Notation.}

\noindent We will use $\Z$ to denote the set of integers and $\Q$ to
denote the set of rational numbers. Also $\Z_D$ will denote the
congruence group modulo $D$.

\noindent For sequences $A$ and $B$, define $A+B :=\{ a+ b | a \in
A, b \in B \}$.

\noindent For an element $b \in \BZ_p$ and a sequence $A$, define $b
\cdot A := \{ba | a \in A \}$.

\noindent A good way to present a sequence  $A$  is to write
$A:=\{a_1^{[m_1]}, \dots, a_k ^{[m_{k}]} \},$ where $m_a$ is the
multiplicity of $a$ in $A$ (sometime we use the notation $m_a(A)$ to
emphasize the role of $A$), and $a_1, \dots, a_k$ are the different
elements of $A$.

\noindent The maximum multiplicity of $A$ is $m(A):=\max_{a \in
\BZ_p} m_a(A).$ We will always assume that $m(A) \le p$, for every
sequence $A$ in the paper.

\noindent We say $A$ is decomposed into subsequences $A_1,\dots,A_k$
and write $A=\bigcup_{i=1}^{*k} A_i$ if $m_a(A)=\sum_{i=1}^k
m_a(A_i)$ for every $a\in \Z_p$.

\vskip2mm

\noindent Asymptotic notation will be used under the assumption that
$p \rightarrow \infty$. For $x \in \BZ_p$, $\|x\|$ (the norm of $x$)
is the distance from $x$ to $0$. (For example, the norm of $p-1$ is
1).

\noindent A subset $X$ of $\Z_p$  is called a  $K$-net if for any
$n\in \Z_p$ there exists $x\in X$ such that $n\in [x,x+K]$. It is
clear that if $X$ is a $K$-net, then $X+T=\Z_p$ for any interval $T$
of length $K$ in $\Z_p$. We will use the same notion over $\Z$ and
$\Q$ as well.

\noindent For a finite set $X$ of real numbers we use $\min(X)$(or,
respectively, $\max(X)$) to denote the minimum (respectively,
maximum) element of $X$.

\section {The classifications.}\label{section:result}

In order to make the statements of the  theorems less technical, we
define

$$f(p,m):=\left \lfloor (pm)^{6/13}\log^2 p  \right \rfloor.$$

\subsection {Zero-sum-free sequences.} View the elements of
$\BZ_p$ as integers between $0$ and $p-1$. The most natural way to
construct a zero-sum-free sequence is to select non-zero elements
whose sum is less than $p$. Our first theorem shows that this is
essentially the only way.

\begin{theorem}\label{theorem:main1} There is a positive constant
$c_1$ such that the following holds. Let $1\le m\le p$ be a
positive integer and $A$ be a zero-sum-free sequence of $\BZ_{p}$
satisfying $m(A)\le m$. Then there is a non-zero residue $b$ and a
subsequence $A^{\flat}\subset A$ of cardinality at most
$c_1f(p,m)$ such that

$$\sum_{a\in b\cdot (A\backslash A^{\flat})}a<p.$$

\end{theorem}

\noindent Notice that zero-sum-freeness and incompleteness are
preserved under dilation. This explains the presence  of the element
$b$ in the theorem. Another issue one needs to address is the
cardinality of the exceptional sequence $A^{\flat}$. It is known
(and not hard to prove) that most zero-sum-free sequences with
maximum multiplicity $m$ in $\Z_p$ have cardinality $\Theta
((pm)^{1/2})$. Thus, in most cases, the cardinality of $A^{\flat}$
(which is at most $(pm)^{6/13+o(1)}$) is negligible compared to that
of $|A|$. (The same will apply for later results.)  Exceptional
sequences cannot be avoided (see Sections
\ref{section:zerofree},\ref{section:incomplete} and also
\cite{NgSzemVu}).

\noindent By setting $m=1$, we have the following corollary for the
case when $A$ is a set.

\begin{corollary}\label{cor:main1} There is an absolute positive constant
$c_1$ such that the following holds. For any {\it zero-sum-free}
subset $A$ of $\BZ_{p}$ there is a non-zero residue $b$ and a  set
$A^{\flat}\subset A$ of cardinality at most $c_1f(p,1)$ such that

$$\sum_{a\in b\cdot  (A\backslash A^{\flat}) }a<p.$$

\end{corollary}

\subsection{Incomplete sequences.}
The easiest way to construct an incomplete sequence is to select
elements with small norms. Clearly, if $A$ is a sequence where
$\sum_{a \in A} \|a \| < p-1$ then $A$ is {\it incomplete}. Our
second theorem shows that this trivial construction is essentially
the only possibility.

\begin{theorem}\label{theorem:main2} There is a positive constant $c_2$
such that the following holds. Let $1\le m\le p$ be a positive
integer and $A$ be an incomplete sequence in $\BZ_{p}$ satisfying
$m(A)\le m$. Then there is a non-zero element $b \in \BZ_p$ and a
subsequence $A^{\flat}\subset A$ of cardinality at most
$c_2f(p,m)$ such that

$$\sum_{a\in b \cdot (A\backslash A^{\flat})}\|a\|<p.$$

\end{theorem}

\noindent By setting $m=1$, we have

\begin{corollary}\label{cor:main2} There is a positive constant $c_2$ such that the following holds.
For any incomplete subset $A$ of $\BZ_p$ there is a non-zero
residue $b$ and a  set $A^{\flat}\subset A$ of cardinality at most
$c_2f(p,1)$ such that

$$\sum_{a\in b \cdot (A\backslash A^{\flat})}\|a\|<p.$$
\end{corollary}

\subsection{$l$-incomplete sequences.}
View $A$ as a sequence of integers in the interval
$[-(p-1)/2,(p-1)/2]$. Our classification in this subsection is a
little bit different from the previous two. We are going to classify
the structure of $\sum_{l}(A)$ instead of that of $A$. The reason is
that this classification is natural and easy to state. Furthermore,
it is also easy to derive information about $A$ using the
classification of $\sum_{l}(A)$.

\noindent If  all $l$-sums of $A$ belong to an interval of length
less than $p$ in $\Z$, then $A$ is $l$-incomplete in $\BZ_p$. Of
course, the converse is not true. However, our third theorem says
that the reversed statement can be obtained at the cost of a small
modification (in the spirit of the previous theorems).

\begin{theorem}\label{theorem:main3} There is a positive constant
$c_3$ such that the following holds. Let $1\le m\le p$ be a positive
integer, let $A$ be a sequence in $\Z_p$, and let $l$ be an integer
satisfying  $c_3f(p,m) \le l \le |A|-c_3f(p,m).$ Assume furthermore
that $A$ is $l$-incomplete and $m(A)\le m$. Then there exist

\begin{itemize}
\item residues $b,c \in \BZ_p$ with $b\neq 0$,
\item a sequence $A^{\flat}\subset A$ of cardinality less than $c_3f(p,m),$ and
\item an integer $l_1\ge l-2f(p,m)$
\end{itemize}

\noindent such that the union $\bigcup_{l_1\le l' \le
l_1+(pm)^{3/13}}\sum_{l'}(A')$ is contained in an interval of length
less than $p$, where $A':=b \cdot (A\backslash A^{\flat}) + c$ is
considered as a sequence of integers in $[-(p-1)/2,(p-1)/2]$.
\end{theorem}

\noindent The property $l$-incompleteness is preserved under linear
transforms. This explains why we need two parameters $b$ and $c$ in
the theorem. The reader is invited to state a corollary for the case
when $A$ is a set.

\section{Structure of long zero-sum-free
sequences.}\label{section:zerofree}

\noindent Let $1\le m \le p$ be a positive integer and  $A$ be a
zero-sum-free sequence of $\Z_p$ with maximum multiplicity $m(A) \le
m$. Trying to make $A$ as long as possible, we come up with the
following natural candidate

$$A_{1}^{m} :=\{1^{[m]},2^{[m]},\dots,(n-1)^{[m]}, n ^{[k]}\}$$

\noindent where $k$ and $n$ are the unique integers satisfying $1
\le k \le   m$ and

$$ m(1+2+\dots+n-1)+kn <  p \le   m(1+2+\dots+n-1)+(k+1)n. $$

\noindent As a consequence of Theorem \ref{theorem:main1}, one can
show that any zero-sum free sequence with $m(A) \le m$ and
cardinality close to $|A^m_{1}|$ is almost a subsequence of $A^m_1$,
after a proper dilation.

\begin{theorem}\label{theorem:zerofree:1}
Let $6/13< \alpha < 1/2$ be a fixed constant. Assume that $A$ is a
zero-sum-free sequence of $\Z_p$ with maximum multiplicity $m(A)\le
m$ and cardinality $|A_1^{m}|-O((pm)^{\alpha})$. Then there is a
non-zero element  $b \in \BZ_p$ and a subsequence $A^{\flat}\subset
A$ of cardinality $O((pm)^{(\alpha+1/2)/2})$ such that $b \cdot
(A\backslash A^{\flat})\subset A^m_1$.

\end{theorem}

\noindent We can go further by showing not only that $|A\backslash
A_1^m|$ is small, but also that the sum of the norm of the elements
in this sequence is small. An example is given by Theorem 1.9 of
\cite{NgSzemVu}, which we restate below.

\begin{theorem}\label{theorem:zerofree:2} \cite{NgSzemVu} Let $A$ be a
zero-sum-free subset of $\BZ_p$  of size at least $.99 \sqrt{2p}$.
Then there is some non-zero element $b \in \BZ_p$ such that
$$\sum_{a \in b \cdot A, a<p/2} \| a\| \le p+O(p^{1/2}) $$
and
$$\sum_{a \in b \cdot A, a>p/2} \| a\|=O(p^{1/2}).$$
\end{theorem}

\noindent The bound $O(p^{1/2})$ is sharp (see \cite{D1} for further
discussion).

\noindent Now assume that  the cardinality of $A$ differs from that
of the extreme example $A^m_1$ by a constant. In this case, we can
tell exactly what $A$ is.

\noindent Let $n(p)$ denote the largest integer $n$ such that
  $$\sum_{i=1}^{n-1}i < p.$$

\begin{theorem} \label{theorem:zerofree:3} \cite{NgSzemVu} There is a constant
$C$ such that the following holds for all primes $p \ge C$.

\begin{itemize}

\item  If $p \neq \frac{n(p)(n(p) +1)}{2} -1$, and $A$ is a subset
of $\BZ_p$ with $n(p)$ elements, then $0 \in \sum(A)$.

\item If $p = \frac{n(p)(n(p) +1)}{2} -1$, and $A$ is a subset of
$\BZ_p$ with $n(p)+1$ elements, then $0 \in \sum(A)$. Furthermore,
up to a dilation, the only zero-sum-free set with $n(p)$ elements is
$\{-2,1,3,4, \dots, n(p)\}$.

\end{itemize}

\end{theorem}

{\it Remarks.} Theorem \ref{theorem:zerofree:3} is also obtained
independently by J. M. Deshouillers and G. Prakash.

\noindent We sketch the proof of  Theorem \ref{theorem:zerofree:1}.

\begin{proof} (Proof of Theorem \ref{theorem:zerofree:1}.) Theorem \ref{theorem:main1}
implies that there is a non-zero residue $b$ and a subsequence
$A^{\flat}\subset A$ of cardinality less than $c_1f(p,m)$ such that
$\sum_{a \in A'} a<p,$ where $A'=b \cdot (A \backslash A^{\flat})$
is viewed as sequence of integers in $[1,p-1]$.

\noindent Notice that $|A'|=|A_1^{m}|-O((pm)^{\alpha})-c_1f(p,m)=
|A^1_m|-O((pm)^\alpha)$. For short put $t=|A' \backslash A_1^m|.$ It
follows from the inequality $n + \sum_{a\in A_1^m} a \ge p \ge
\sum_{a \in A'} a$ that

\begin{equation}\label{equation:application:zero:1}
\sum_{a \in A' \backslash A_1^m} a \le n+ \sum_{a \in A_1^m
\backslash A'} a.
\end{equation}

\noindent Let $A_1'$ be the any subsequence of cardinality $t$ in
$A_1^m \backslash A'$ and let $A_1''=A_1^m \backslash (A'\cup
A_1')$. Note that

$$|A_1''|=|A_1^m|-|A'|= O(pm)^{\alpha} \mbox{ and } a\le n\le (2p/m)^{1/2}+1$$

\noindent for any $a\in A_1''$. Thus

\begin{equation}\label{equation:application:zero:2}
n + \sum_{a\in A_1''} a = O(pm)^{\alpha}(p/m)^{1/2}.
\end{equation}

\noindent On the other hand, by definition, every element of $A'
\backslash A_1^m$ is strictly greater than every element of $A_1'$.
Additionally, since the maximum multiplicity is $m$, we have

$$\sum_{a \in A' \backslash A_1^m} a - \sum_{a \in A_1'} a \ge 1+\cdots+1+2+\dots+2+3+\cdots+3+\cdots,$$

\noindent where on the right hand side all numbers (with the
possible exception of the last) appear exactly $m$ times and the
total number of summands is $t$. It is clear that such a sum is
greater than $t^2/3m$; thus

\begin{equation}\label{equation:application:zero:3}
\sum_{a \in A' \backslash A_1^m} a - \sum_{a \in A_1'} a\ge t^2/3m.
\end{equation}

(\ref{equation:application:zero:1}),(\ref{equation:application:zero:2}),(\ref{equation:application:zero:3})
together give

$$t^2/3m \le \sum_{a \in A' \backslash A_1^m} a - \sum_{a \in A_1'} a
\le n + \sum_{a\in A_1''} a = O(pm)^{\alpha}(p/m)^{1/2}.$$

\noindent In other words, $t=O((pm)^{(\alpha+1/2)/2})$.

\end{proof}

{\it Remarks.} The interested reader may also read \cite[Section
7]{Ge} and \cite{GLP1,GLP2} for further results on long
zero-sum-free sequences.

\section{Structure of long incomplete sequence.}\label{section:incomplete}

\noindent Let $1\le m \le p$ be a positive integer and  $A$ be an
incomplete sequence of $\Z_p$ with maximum multiplicity $m(A)\le m$.
Trying to make $A$ as large as possible, we come up with the
following example,

$$A_2^{m}=\{-n^{[k]}, -(n-1)^{[m]}, \dots,-1^{[m]},0^{[m]},1^{[m]},\dots,(n-1)^{[m]}, n^{[k]}\}$$

\noindent where $1 \le k \le   m$ and $n$ are the unique integers
satisfying

$$2m(1+2+\dots+n-1)+2kn<p\le 2m(1+2+\dots+n-1)+2(k+1)n. $$

\noindent Using Theorem \ref{theorem:main2}, we can prove the
following.

\begin{theorem}\label{theorem:incomplete:1} Let $6/13< \alpha< 1/2$
be a fixed constant. Assume that $A$ is an incomplete sequence of
$\Z_p$ with maximum multiplicity $m$ and cardinality $|A|=
|A_2^{m}|-O((pm)^{\alpha})$. Then there is a non-zero element $b
\in \BZ_p$ and a subsequence $A^{\flat}\subset A$ of cardinality
$O((pm)^{(\alpha+1/2)/2})$ such that $b \cdot (A\backslash
A^{\flat})\subset A^m_2$.
\end{theorem}

\noindent The proof is similar to that of Theorem
\ref{theorem:zerofree:1} and is omitted.

\noindent As an analogue of Theorem \ref{theorem:zerofree:2}, we
have

\begin{theorem} \label{theorem:incomplete:2} \cite{NgSzemVu}  Let $A$ be an incomplete
subset of $\BZ_p$  of size at least $1.99 p^{1/2}$. Then there is
some non-zero element $b \in \BZ_p$ such that

$$\sum_{a \in b \cdot A} \| a\| \le p +O(p^{1/2}). $$

\end{theorem}

(Again, the error term $O(p^{1/2})$ is sharp, see \cite{D2} and
\cite{DF}.)

\noindent A well-known theorem of {J. E. Olson} \cite{Olson2} gives
a sharp estimate for the maximum cardinality of an incomplete set.

\begin{theorem}\label{theorem:incomplete:3}
Let $A$ be a subset of $\Z_p$ of cardinality more than
$(4p-3)^{1/2}$. Then $A$ is complete.
\end{theorem}

\section{The number of zero-sum-free and incomplete
sequences.}\label{section:partition}

 \noindent In this section we apply Theorems \ref{theorem:main1}, \ref{theorem:main2} to count the number of
 zero-sum-free sequences and incomplete sequences.

 \noindent We fix $m$. The following theorem is well known in theory of partitions (a corollary of a theorem of G. Meinardus, \cite[Theorem 6.2]{A}).

 \begin{theorem}\label{theorem:partition:m} Let $p_m(n)$ be the number of partitions of $n$ in which each
 positive integer appears at most $m$-times. Then

 $$p_m(n)=\exp{((\sqrt{(1-\frac{1}{m+1})\frac{2}{3}}\pi + o(1))\sqrt{n})}.$$

 \end{theorem}

 \noindent By Theorem \ref{theorem:main1}, the main part of zero-sum-free sequences (after a proper dilation)
 corresponds to a partition of a number less than $p$. Thus, using Theorem
 \ref{theorem:partition:m}, we infer the following.

 \begin{theorem}\label{theorem:zerofree:partition:m}
 Let $N_1^m$ be the number of zero-sum-free sequences $A$ satisfying $m(A) \le
 m$. Then
 $$N_1^m = \exp((\sqrt{(1-\frac{1}{m+1})\frac{2}{3}}\pi + o(1))\sqrt{p}).$$
 \end{theorem}

 \begin{corollary}\label{cor:zerofree:partition:set}
 The number of zero-sum-free sets is
 $\exp((\sqrt{\frac{1}{3}}\pi + o(1))\sqrt{p})$.
 \end{corollary}

 \noindent By Theorem \ref{theorem:main2}, the main part of incomplete
 sequences (after a proper dilation) can be split into two parts, each of which corresponds to
 a partition of a number less than $p$. Thus we obtain the following.

 \begin{theorem}\label{theorem:incomplete:partition:m}
 Let $N_2^m$ be the number of incomplete sequences $A$ satisfying $m(A) \le
 m$. Then
 $$N_2^m = \exp((\sqrt{(1-\frac{1}{m+1})\frac{4}{3}}\pi + o(1))\sqrt{p}).$$
\end{theorem}

 \begin{corollary}\label{cor:incomplet:partition:set}
 The number of incomplete sets is
 $\exp((\sqrt{\frac{2}{3}}\pi + o(1))\sqrt{p})$.
 \end{corollary}

 \begin{proof}(Proof of Theorem \ref{theorem:zerofree:partition:m})
 The lower bound for $N_1^m$ is obvious, any partition of $p-1$ in which each number appears at most
 $m$-times gives a zero-sum-free sequence of maximum multiplicity bounded by $m$.

 \noindent For the upper bound, we apply Theorem \ref{theorem:main1}. First, the
 number of choice for $A^\flat$ is $\sum_{n\le (pm)^{6/13+o(1)}}\binom{pm}{n}=\exp({o(\sqrt
 p)})$. Second, the elements of $A':=b(A\backslash A^\flat)$ forms a partition of
 $\sum_{a\in A'} a$ (which is a positive integer less than $p$) in which each positive integer appears at most $m$-times.
 Hence, the number of choice for $A^\sharp$
 is at most

 $$\sum_{n\le p-1} p_m(n) \le p \exp((\sqrt{(1-\frac{1}{m+1})\frac{2}{3}}\pi +
 o(1))\sqrt{p}).$$

 \noindent Finally, together with dilations, the number of
 zero-sum-free sequences is bounded by

 $$p^2\exp((\sqrt{(1-\frac{1}{m+1})\frac{2}{3}}\pi + o(1))\sqrt{p}) = \exp((\sqrt{(1-\frac{1}{m+1})\frac{2}{3}}\pi + o(1))\sqrt{p}).$$

 \end{proof}

 \begin{proof}(Proof of Theorem \ref{theorem:incomplete:partition:m})
 The lower bound for $N_2^m$ is again obvious, any two partitions of $(p-3)/2$
 in which each number appears at most $m$-times give two
 nonnegative sequences. We then take the union of one sequence with
 the negative of the other sequence. It is not hard to check that
 the formed sequence $A$ is incomplete and $m(A) \le m$. Thus
 $$N_2^m \ge (p_m((p-1)/2))^2 = \exp((\sqrt{(1-\frac{1}{m+1})\frac{4}{3}}\pi + o(1))\sqrt{p}).$$

 \noindent For the upper bound we use Theorem \ref{theorem:main2}. Argue similarly as in the proof of Theorem \ref{theorem:zerofree:partition:m},
 we infer that the number of exceptional sequences $A^\flat$ is at most
 $e^{o(\sqrt{p})}$. Write $A':=b(A\backslash A^\flat)=A^{+}\cup A^{-}$, the decomposition
 of $A'$ into sequences of
 nonnegative and negative elements respectively. The elements of $A^+$ form a
 partition of $\sum_{a\in A^+} a$ in which each positive integer appears at
 most $m$-times. The elements of $A^{-}$ corresponds to a partition of $\sum_{a\in A^-}(-a)$ in which each (negative) number appears at most $m$-times.
 Thus the number of choice for $A'$ is at most

 $$\sum_{k+l < p} p_m(k)p_m(l)
 \le p^2 \exp((\sqrt{(1-\frac{1}{m+1})\frac{4}{3}}\pi +
 o(1))\sqrt{p}).$$

 \noindent Putting everything together, we obtain an upper bound for
 $N_2^m$,

 \begin{eqnarray*}
 N_2^m & \le & p e^{o(\sqrt{p})} p^2 \exp((\sqrt{(1-\frac{1}{m+1})\frac{4}{3}}\pi + o(1))\sqrt{p})\\
 & \le & \exp((\sqrt{(1-\frac{1}{m+1})\frac{4}{3}}\pi + o(1))\sqrt{p}).
 \end{eqnarray*}

 \end{proof}

\section{$l$-incomplete sequences}\label{section:l-incomplete}

\noindent Assume that $A,l,m$ satisfy conditions of Theorem
\ref{theorem:main3}. Trying to make $A$ as large as possible, we
come up with the following example,

$$A_3^{m}=\{-n^{[k]}, -(n-1)^{[m]}, \dots,-1^{[m]},0^{[m]},1^{[m]},\dots,(n-1)^{[m]}, n^{[k]}\}$$

\noindent where $k$ and $n$ are the optimal integers such that  $1
\le k \le m$ and all the $l$-sums of $A_3^m$ are contained in an
interval of length less than $p$.

\noindent However, the extremal example for $l$-incomplete
sequences, in general, is not unique (for instance if $l=m=p$ then
any sequence $\{-1^{[n]},0^{[p]},1^{[p-2-n]}\}$ is $l-$incomplete
and of maximum cardinality). Nevertheless, Theorem
\ref{theorem:main3} still allows us to conclude that any
$l-$incomplete sequence of size close to $|A_3^m|$ can be dilated
and translated into one of the extremal examples, as in the spirit
of Theorems \ref{theorem:zerofree:1} and \ref{theorem:incomplete:1}.

 \noindent Let us discuss in detail the special case $l=p$.
 This is motivated by the classical  theorem of  {P. Erd\H{o}s, A. Ginzburg} and {A.
Ziv} \cite{EGZ},  one of the starting points of combinatorial
number theory.

\begin{theorem}\label{theorem:l-zerofree:1} (Erd\H{o}s-Ginzburg-Ziv) For any sequence $A\in Z_p$ of cardinality $2p-1$ there is a
subsequence $A' \subset A$ of cardinality $p$ such that $\sum_{a\in
A'} a = 0$.

\end{theorem}

\noindent In fact, {P. Erd\H{o}s, A. Ginzburg} and {A. Ziv} proved
the statement for any finite abelian group $G$, by reducing it to
the case $G=\Z_p$ above.

\noindent In the context of this paper, Theorem
\ref{theorem:l-zerofree:1} stated that any sequence of cardinality
$2p-1$ in $\Z_p$ is not $p$-zero-sum-free. The bound $2p-1$ is sharp
as shown by the example $A=\{a^{[p-1]},b^{[p-1]} \}$, for any two
different elements $a,b \in \Z_p$. Using Theorem
\ref{theorem:main3}, we prove that  if $A$ is $p$-zero-sum-free and
$|A|-p \gg f(p,p)= \lfloor p^{12/13} \log^2 p \rfloor$, then $A$ has
two elements of high multiplicities.

\begin{theorem}\label{theorem:l-zerofree:2}
There is a positive constant $C$ such that the following holds for
all primes $p>C$. Assume that $A$ is a $p$-zero-sum-free sequence
and $p+c_3f(p,p) \le |A| \le 2p-2$. Then
$\{a^{[m_a]},b^{[m_b]}\}\subset A$, where $a,b$ are two different
elements of $\Z_p$ and $m_a+m_b \ge 2(|A|-p-(c_3+3)f(p,p))$.

\end{theorem}

\noindent Notice that $A$ must have at least $p$ elements so that
the notion of $p$-zero-sum-free makes sense. Our theorem already
yields a non-trivial conclusion  when $A$ has slightly more than $p$
elements. A similar statement was proved in \cite{GPT} (see also
\cite{CS}), but under the stronger assumption that $|A| \ge
\frac{3}{2} p$.

\noindent As a quick application of Theorem
\ref{theorem:l-zerofree:2}, one obtains the following refinement of
Theorem \ref{theorem:l-zerofree:1}, which was first proved by B.
Peterson and T. Yuster.

\begin{corollary}\cite[Section 7]{Ge} \label{cor:EGZ1} The following holds for all sufficiently large
primes $p$.  Let $A$ be a $p$-zero-sum-free sequence of cardinality
$2p-2$ in $\Z_p$. Then $A= \{a^{[p-1]},b^{[p-1]} \}$, where $a,b$
are two different elements of $\Z_p$.
\end{corollary}

\begin{proof} (Proof of Corollary \ref{cor:EGZ1}) By Theorem
\ref{theorem:l-zerofree:2}, we may assume that $$A=\{0^{[p-k_1]},
1^{[p-k_2]}, a_1, \dots, a_l \}$$ where $1\le k_1 =o(p) , 1\le
k_2=o(p), l=k_1+k_2-2$ and $a_i$ are (not necessarily distinct)
integers in $[-p/2, p/2]\backslash \{0,1\}$. If $l=0$ then we are
done. Assume that  $l \ge1$. We are going to construct a subsequence
of $A$ of length $p$ whose elements sum up to zero modulo $p$.

{\bf Case 1:} There is some $a_i$ with absolute value at least
$p/6$.

\noindent Assume that $p/2> a_1 \ge p/6$. The subsequence
$\{0^{[a_1-1]}, 1^{[p-a_1]}, a_1\}$ has cardinality $p$ and sums up
to zero modulo $p$. In the case $-p/2 < a_1 \le -p/6$, consider the
subsequence $\{0^{[p-|a_1|-1]}, 1^{[|a_1|]}, a_1\}$.

{\bf Case 2:} All $a_i$ have absolute value less than $p/6$ and
there are at least $\max \{1, k_1-1\}$ negatives among them.

\noindent By a greedy algorithm, one can find a non-empty sequence
(say, $a_1, \dots, a_{l_1}$) of negative elements  such that $l_1 +
|a_1 + \dots + a_{l_1}| \ge k_1$.  Then the subsequence

$$\{0^{[p- l_1 - |a_1+\dots a_{l_1}|]}, 1^{[|a_1+ \dots a_{l_1}|]}, a_1, \dots, a_{l_1}\}$$

sums up to zero modulo $p$.

{\bf Case 3:} All $a_i$ have absolute value less than $p/6$ and
there are at least $\min \{l, k_2\}$ positives among them.

\noindent As each positive element is at least 2 and at most $p/6$,
there is a subsequence of (say, $l_2$) positive elements whose sum
is at least $k_2$ and at most $p/3$. Assume that $a_1, \dots,
a_{l_2}$ are these elements. Then the subsequence $$\{0^{[(a_1+
\dots a_{l_1})-l_2]}, 1^{[p- (a_1+ \dots a_{l_2})]}, a_1, \dots,
a_{l_2} \}$$ sums up to zero modulo $p$.
\end{proof}

\noindent We conclude this section by sketching the proof of Theorem
\ref{theorem:l-zerofree:2}.

\begin{proof} (Sketch of proof of Theorem \ref{theorem:l-zerofree:2})
Since $A$ is $p$-zero-sum-free in $\Z_p$, $A$ is also
$p$-incomplete. By Theorem  \ref{theorem:main3}, after a linear
transform, we can find a subsequence $A'$ of $A$ such that

\begin{equation}\label{equation:application:1}
\max \{\sum_{l_1}(A')\}-\min \{\sum_{l_1}(A')\}<p,
\end{equation}

\noindent where $l_1\ge p-2f(p,p)$ and $|A'|\ge |A|-c_3f(p,p)$ and
where $c_3$ is a positive constant.  (Recall that $\max (X) $
(respectively, $\min (X)$) refers to the maximum (respectively,
minimum) element in $X$.)

\noindent Let $A'=\{a_1,\dots, a_q\}$, where $a_i\le a_{i+1}$ for
$1\le i \le q-1=|A'|-1$ and rewrite (\ref{equation:application:1})
as

\begin{equation}\label{equation:application:2}
\sum_{i=1}^{l_1} a_{q-l_1+i}- \sum_{i=1}^{l_1} a_i = \sum_{i=1}^k
a_{q-k+i}- \sum_{i=1}^k a_i < p,
\end{equation}

\noindent where $k=\min (l_1,q-l_1)$. Note that

\begin{equation}\label{equation:application:3}
\sum_{i=1}^k a_{q-k+i}- \sum_{i=1}^k a_i \ge \sum_{i=i_0}^{j_0}
a_{i+p} - \sum_{i=i_0}^{j_0}a_i=\sum_{i=i_0}^{j_0} (a_{i+p}-a_i),
\end{equation}

\noindent where $i_0=\max(1,q-l_1-p+1)$ and $j_0=\min(l_1,q-p)$.

\noindent Since $A$ has maximum multiplicity less than $p$, we have,
for any $i$, that $a_{i+p}-a_i \ge 1. $ Thus by
(\ref{equation:application:3}) we obtain that

$$j_0-i_0 \le \sum_{i=i_0}^{j_0} (a_{i+p}-a_i) < p,$$

\noindent and we infer that the number of $i\in [i_0,j_0]$ such that
$a_{i+p}-a_i=1$ is at least  $2(j_0-i_0)-p+3$. Next let $i_1$ and
$j_1$ be the smallest and largest index $i$ in $[i_0,j_0]$ such that
$a_{i+p}-a_i=1$. Thus $a_{i_1+p}-a_{i_1}=a_{j_1+p}-a_{j_1}=1$ and

\begin{equation}\label{equation:application:4}
2(j_0-i_0)-p+2 \le j_1-i_1 \le j_0-i_0 <p.
\end{equation}

\noindent In what follows, $a_{i_1}$ plays a special role, so we
denote it by $a$ to distinguish it from the other $a_i$. Let
$B=\{a_{i_1},\dots,a_{j_1+p}\}.$ Obviously $|B|=j_1-i_1+p+1$ and
$a_{j_1+p}-a_{i_1} \le 2$.

\noindent Set $\gamma := j_0-i_0$. Then $0\le \gamma \le l_1-1$. We
consider two cases.

{\bf Case 1: $a_{j_1}=a$.} In this case $a_{j_1+p}=1$ and
$B=\{x^{[m_0]},(x+1)^{[m_1]}\}$ where

\begin{equation}\label{equation:application:5}
m_0+m_1=j_1-i_1+p+1 \ge 2(j_0-i_0)-p+2+p+1=2\gamma+3.
\end{equation}

{\bf Case 2: $a_{j_1}=a+1$.} Recall that the number of pairs
$(a_i,a_{i+p})$ such that $a_{i+p}-a_i=1$ is at least $2(j_0-i_0)-p
+2= 2\gamma - p+2$. Furthermore if $a_{i+p}-a_i=1$ then either $a_i$
or $a_{i+p}$ must be $a+1$. By this observation, none of the
elements in $\{a_{j_1+1},\dots, a_{p+i_1-1}\}$ belongs to any pair
$(a_i,a_{i+p})$ with $a_{i+p}-a_i=1$. Furthermore, we have $a_i=a+1$
for $j_1+1\le i \le p+i_1-1$. As a consequence, the multiplicity
$m_1$ of $a+1$ in $B$ is at least

\begin{equation}\label{equation:application:6}
m_1 \ge 2\gamma-p+2+(p+i_1-j_1-1) = 2\gamma -(j_1-i_1)+1.
\end{equation}

\noindent It is convenient to write $B=
\{a^{[m_0]},(a+1)^{[m_1]},(a+2)^{[m_2]}\}$. Clearly we have $\min
(p^\ast B)= \min(p-m_0,m_1) + 2 (p-m_0-\min(p-m_0,m_1))$ and $\max
(p^\ast B) = 2m_2 + \min (p-m_2,m_1)$.

\noindent Besides, it is not hard to show that

\begin{equation}\label{equation:application:7}
\sum_p(B)=[\min (\sum_p(B)),\max (\sum_p(B) )].
\end{equation}

\noindent The $p$-zero-sum-free assumption implies that $\max
(\sum_p(B)) < p$. It follows that
\begin{equation}\label{equation:application:8}
2m_2 + \min (p-m_2,m_1) < p.
\end{equation}

\noindent Consequently,

\begin{equation}\label{equation:application:9}
2m_2 + m_1 < p.
\end{equation}

\noindent From (\ref{equation:application:6}) and
(\ref{equation:application:9}) we deduce that $m_2 \le (p- 2\gamma +
(j_1-i_1)-2)/2.$ On the other hand, $m_0+m_1+m_2 = |B|=
j_1-i_1+p+1$. Thus

$$m_0+m_1 \ge j_1-i_1+p+1 - (p- 2\gamma + (j_1-i_1)-2)/2 \ge \gamma
+2+(j_1-i_1+p)/2.$$

\noindent The latter inequality, together with
(\ref{equation:application:4}), yields

\begin{equation}\label{equation:application:10}
m_0+m_1 \ge 2 \gamma +3.
\end{equation}

\noindent To summarize, in both cases
((\ref{equation:application:5}) and (\ref{equation:application:10}))
we have $m_0+m_1 \ge 2\gamma + 3$. Combining this with the estimates
$l_1 \ge p- 2f(p,p)$ and $q \ge |A|- c_3f(p,p)$ we get

\begin{eqnarray*}
m_0+m_1 & \ge & 2 (\min(l_1,q-p)-\max(1,q-l_1-p+1)) +3 \\
& \ge & 2(|A|-p)-(2c_3+6) f(p,p).
\end{eqnarray*}
\end{proof}

\section{The key lemmas.}\label{section:lemma}

\noindent The key lemmas we use in proofs are the following results
from \cite{SzemVu2}.

\begin{theorem}\label{lemma:main1} For any fixed positive integer $d$ there exist positive $C=C(d)$ and
$c=c(d)$ depending on $d$ such that the following holds. If $A$ is a
subset of $[n]$ and $l$ is a positive integer such that $l^d|A|\ge
C(d)n$ and $l\le |A|/2$. Then $\sum_{l}(A)$ contains an arithmetic
progression of length $c(d)l|A|^{1/d}$.
\end{theorem}

\vskip2mm

\begin{theorem}\label{lemma:main2} For any fixed positive integer $d$ there exist positive $C=C(d)$ and
$c=c(d)$ depending on $d$ such that the following holds. If $A$ is a
subset of $\BZ_p, |A|\ge 2$ and $l$ is a positive integer such that
$l^{d+1}|A|\ge C(d)p$, then $\sum_{l}(A)$ contains all residue
classes modulo $p$ or contains an arithmetic progression of length
$c(d)l|A|^{1/d}.$
\end{theorem}

\vskip2mm

\begin{theorem}\label{lemma:main3} For any fixed positive integer $d$ there exist positive $C=C(d)$ and
$c=c(d)$ depending on $d$ such that the following holds. Let
$A_1,\dots,A_l$ be subsets of cardinality $|A|$ of $\BZ_p$ where $l$
and $|A|$ satisfy $l^{d+1}|A|\ge C(d)p$. Then $A_1+\dots+A_l$
contains all residue classes modulo $p$ or an arithmetic progression
of length $c(d)l|A|^{1/d}$.
\end{theorem}

\noindent In our proofs we  will be mainly interested in the case
$d=1$ and $d=2$. We will also use the  following lemmas. The proofs
are left as exercises.

\begin{lemma}\label{lemma:simple1} \cite{Olson2}
There are positive constants $C_0$ and $c_0$ such that the following
holds. Let $A$ be a set of $\Z_p$ satisfying $|A| \le C_0p^{1/2}$.
Then
$$|\sum_{l}(A)|\ge c_0|A|^2$$
where $l=\lfloor |A|/2 \rfloor$.
\end{lemma}

\vskip2mm

\begin{lemma}\label{lemma:simple2}
Let $D$ be a positive integer and $X$ be a sequence of cardinality
$D$ in $\Z_D$. Then $\sum(X)$ contains the zero element.
Furthermore, if the elements of $X$ are co-prime with $D$, then
$\sum(X) = \Z_D$.
\end{lemma}

\vskip2mm

\begin{lemma}\label{lemma:simple3}\cite{SzemVu1}
Let $d_1,\dots,d_n$ be distinct positive integers and
$D=lcm(d_1,\dots,d_n)$. Then for any $0 \le r \le D-1$ there exist
$0\le a_i\le d_i-1$ such that $\sum_{i=1}^n a_i/d_i=r/D (\mod \ 1)$.
\end{lemma}

\begin{lemma}\label{lemma:simple4} (a consequence of
Chinese remainder theorem) Let $d_1,\dots,d_n,D$ be distinct
positive integers and $\gcd(d_1,\dots,d_n,D)=1$. Then for any $0 \le
r \le D-1$ there exist $0\le a_i\le D$ such that $\sum_{i=1}^n a_i
\le D$ and $\sum_{i=1}^n a_id_i/D=r/D(\mod 1)$.
\end{lemma}

\vskip2mm

\noindent We  will mainly focus on the proof of Theorem
\ref{theorem:main3}, which is the most difficult among the three
theorems in Section \ref{section:result}. Theorem
\ref{theorem:main2} can be proved by invoking the same technique in
a simpler manner and we will sketch its proof. Theorem
\ref{theorem:main1} can be deduced from Theorem \ref{theorem:main2}
by several applications of Lemma \ref{lemma:main1}.

\section{ Proof of Theorem \ref{theorem:main3}}\label{section:main3}

\noindent Our plan consists of four main steps

\begin{itemize}

\item We first obtain a long arithmetic progression (say $P$) by
using the subset sums of a small subsequence of $A$.

\item Next we show that (after a linear transform)
one can find a reasonably short interval (say $A_0$) around 0 which
contains many elements of $A$.

\item Since $A$ is $l$-incomplete, the sum of the subset sums of
the remaining part $A\backslash (A_0\cup P)$ with $A_0$ and $P$ does
not cover $\Z_p$. Thus the main part of $A$ concentrates around a
few points which are evenly distributed in $\Z_p$.

\item Finally we use this structural information to deduce
the statement of the theorem.

\end{itemize}

\subsection {Creating a long arithmetic progression}

\noindent Assume that $A$ is an $l$-incomplete sequence with maximal
multiplicity less than $m$. Recall that

$$f(p,m)=\lfloor (pm)^{6/13}\log^2 p \rfloor .$$

\noindent In what follows, we think of $m$ and $p$ as fixed and use
shorthand $f$ for $f(p,m)$. By setting $c_3$ large, we can assume
that  $|A|/f$ is large, whenever needed. If there is an element $a$
such that $m_a(A) \ge |A|-f$ then the theorem is trivial, as we can
take $A^{\flat}=\{b\in A, b\neq a\}$. Thus we can assume that $m(A)
<|A|-f$.

\noindent Let $\lambda$ be a sufficiently large constant. We execute
the first step of the plan by showing the following.

\begin{lemma}\label{lemma:main3:3}
There is a subsequence $A^{\flat}\subset A$ of cardinality at most
$f$ whose $l^{\flat}$-sums, for some integer $l^\flat\le f$, contain
an arithmetic progression of length $\lambda (pm)^{12/13}/m$.
\end{lemma}

\noindent Here we abuse the notation  $A^{\flat}$ slightly. The
current $A^\flat$ is not necessarily the $A^\flat$  in Theorem
\ref{theorem:main3}. However, as the reader will see, the latter
will be the union of  the current $A^\flat$ with a very small
sequence of $A$.

\begin{proof} (Proof of Lemma \ref{lemma:main3:3}) We consider three cases.

{\bf Case 1:} $m>(pm)^{6/13}$.

\noindent Since $m(A) \le |A|-f$ by assumption, we can find in $A$
$f$ disjoint sets $A_1, \dots, A_f$, each has exactly two different
elements. Let  $A'= A \backslash \cup_{i=1}^f A_i$. By the
assumption $m > (pm)^{6/13}$, it follows that for each $i=1,\dots,
f$,

$$f^2|A_i| = 2f^2
> (pm)^{12/13} \gg p. $$

\noindent Thus we can  apply Theorem \ref{lemma:main3} to the $f$
sets $A_1,\dots,A_f$ and conclude that their sum $A_1+\dots+A_f$
contains an arithmetic progression $P$ of length $|P| \ge c(1)f|A_i|
> c(1)(pm)^{6/13}\log^2 p$, for some positive constant $c(1)$.

\noindent On the other hand, the assumption $m > (pm)^{6/13}$ yields
that $(pm)^{6/13} \ge (pm)^{12/13}/m$. Thus

$$|P| \ge \lambda (pm)^{12/13}/m $$

\noindent for any fixed  constant $\lambda$. We complete by letting
$A^{\flat}=\bigcup_{i=1} A_i$ and $l^{\flat}=f$.

{\bf Case 2:} $p^{1/5}<m \le (pm)^{6/13}$.

\noindent Let $A^{\flat}$ be an arbitrary  subsequence of
cardinality $f$ in $A$. Since $m(A^{\flat}) \le m(A) \le m$, we can
find in $A^\flat$ disjoint sets $A_1, \dots, A_m$ each of which has
cardinality

$$\lfloor |A_i|=|A^{\flat}|/m \rfloor = \lfloor f/m \rfloor.$$

\noindent Let $k=\lfloor |A_1|/2 \rfloor$. Since $|A_i|\ll p^{1/2}$,
by Lemma \ref{lemma:simple1} we have

$$|\sum_k(A_i)|\ge c_0|A_i|^2.$$

\noindent Next choose a set $B_i$ of cardinality $|B_i| =
c_0|A_i|^2$ from $\sum_k(A_i)$ for all $i$. Since

$$m^2|B_i|\ge m^2c_0(\frac{f}{m}-1)^2 > c_0m^2 \frac{f^2}{4m^2} > (pm)^{12/13}> p^{12/13+2/13}\gg p,$$

\noindent we can apply Theorem \ref{lemma:main3} to the $m$ sets
$B_1, \dots, B_m$ to conclude that the sumset $B_1+\dots+B_m$
contains an arithmetic progression $P$ of length

$$|P|=c(1)m|B_i| = c(1)c_0 m |A_i|^2
> \frac{c(1)c_0}{4}  m \frac{f}{m}^2 > \frac{\lambda(pm)^{12/13}}{m},$$

\noindent for any fixed $\lambda$, thanks to the definition of
$f=f(p,m)$.

\noindent Let $l^{\flat}=mk$. Note that the arithmetic progression
$P$ is contained in $\sum_k(A_1)+\dots+ \sum_k(A_m)$. But the latter
sumset is  a subset of $\sum_{l^{\flat}}(A^{\flat})$. Thus the set
$\sum_{l^{\flat}}(A^{\flat})$ contains an arithmetic progression $P$
of length $|P| \ge \lambda(pm)^{12/13}/m$.

{\bf Case 3:} $m\le p^{1/6}$.

\noindent Again let  $A^{\flat}$ be an arbitrary  subsequence of
cardinality $f$ of $A$. For each element $a$, let $m_a$ be its
multiplicity in $A^{\flat}$. We partition $A^{\flat}$ according the
magnitudes of these multiplicities.  For $0\le i \le \log m -1 $,
let $n_i$ be the number of element $a$ of $A^\flat$ such that
$2^i\le m_a < 2^{i+1}$. It is easy to see that $f=|A^\flat|\le
\sum_{i=0}^{\log m -1 } n_i 2^{i+1}$ (here the $\log$ has base 2),
which implies that there exists an index $0\le i_0 \le \log m -1 $
satisfying

\begin{equation}\label{equation:main3:ap1}
n_{i_0}2^{i_0+1}\ge \frac{f}{\log m}.
\end{equation}

\noindent Let $a_1,\dots,a_{n_{i_0}}$ be elements of $A^\flat$ whose
multiplicity belongs to $[2^{i_0},2^{i_0+1})$. Set
$B_1:=\dots=B_{2^{i_0}}:=\{a_1,\dots, a_{n_{i_0}}\}.$ Then the union
of the $ B_j$ is a subsequence of $A^\flat$. Furthermore,

\begin{equation}\label{equation:main3:ap2}
|B_1|=n_{i_0} \ge \frac{f}{2^{i_0+1}\log m } > \frac{(pm)^{6/13}}{
m}
\end{equation}

\noindent because $2^{i_0}\le m \le p$. Let $l_1=\lfloor |B_1|/2
\rfloor$. By the assumption $m \le p^{1/6}$ we have

$$l_1^2|B_1| > (pm)^{18/13}/(8m^3) \gg p.$$

\noindent Theorem \ref{lemma:main2} applied to $B_1$ with $d=1$,
yields an arithmetic progression $P_1\subset l_1^\ast B_1$ of length

$$|P_1|\ge c(1)l_1|B_1| >  c(1)|B_1|^2/4.$$

\noindent Since each $B_i$ is a duplicate of $B_1$, we obtain
$2^{i_0}$ duplicates $P_1,P_2,\dots,P_{2^{i_0}}$ of $P_1$ in
$l_1^\ast B_1,\dots, l_1^\ast B_{2^{i_0}}$ respectively. Now
consider $P=P_1+\dots+P_{2^{i_0}}$. Notice that

$$|P|=2^{i_0}|P_1|-(2^{i_0}-1) \ge 2^{i_0}|P_1|/2.$$

\noindent By  (\ref{equation:main3:ap1}) and
(\ref{equation:main3:ap2}), we have

$$|P|\ge 2^{i_0}c(1)|B_1|^2/8 = c(1) 2^{i_0}n_{i_0}|B_1|/8 \ge$$

$$\ge (c(1)/8)(f/(2\log m))((pm)^{6/13}/m) > \lambda(pm)^{12/13}/m$$

\noindent for any fixed $\lambda$. Now observe that

$$P \subset \sum_{l_1}(B_1) + \dots +\sum_{l_1}(B_{2^{i_0}}) \subset
\sum_{2^{i_0}l_1}(A^\flat).$$

\noindent Thus by setting $l^\flat = 2^{i_0}l_1$ we conclude that
the collection of $l^\flat$-sums of $A^\flat$ contains an arithmetic
progression of length $\lambda(pm)^{12/13}/m$.

\end{proof}

\noindent By a dilation of $A$ with some nonzero $b'\in \Z_p$, we
can assume that the arithmetic progression $P$ obtained by Lemma
\ref{lemma:main3:3} is an interval, $P=[p_0,p_0+L]$ for some residue
$p_0$ and $L\ge \lambda (pm)^{12/13}/m$.

\subsection{Dense subsequence around zero}

\noindent Let $Q=\lfloor(pm)^{3/13}\rfloor$ and $A'=A \backslash
A^{\flat}$.

\begin{lemma}\label{lemma:main3:4}
There exists a residue $c'\in \Z_p$ such that $(A'+c')\cap
[-p/(2Q^2),p/(2Q^2)]$ contains a subsequence of cardinality $3Q$.
\end{lemma}

\begin{proof} (Proof of Lemma \ref{lemma:main3:4})
Call a pair $(x,y)$ of $\Z_p\times \Z_p$ {\it nice} if

$$p/Q^2< \|y-x\|<L.$$

\noindent Note that if $(x,y)$ is a nice pair then $x+P \cap y+P
\neq \emptyset$ and  $x+P \cup y+P$ is  an interval of length

\begin{equation}\label{equation:main3:1}
|x+P \cup y+P| \ge \min (|P|+ p/Q^2,p).
\end{equation}

\noindent Assume that $B = \{ x_1,y_1, \dots, x_r,y_r \}$ is a
(maximal) sequence of nice pairs in $A'$ (this means that there is
no more nice pair left in $A' \backslash B$). We are going to show
that $r < Q^2$. Assume otherwise. By (\ref{equation:main3:1}),

$$P'=\bigcup_{z_i\in \{x_i,y_i\}, 1\le i \le Q^2}
z_1+\dots+z_{Q^2}+P=\Z_p.$$

\noindent On the other hand, by the assumption of the Theorem,

$$\left |A'\backslash \bigcup_{i=1}^{Q^2}\{x_i,y_i\} \right | = |A|-|A^\flat|-2Q^2 \ge
|A|-2f  \ge l.$$

\noindent So we are able to choose a subsequence $C$ in
$A'\backslash B$ of cardinality $l-l^\flat-Q^2$.

\noindent But then

$$\Z_p= P'+\sum_{c\in C} c \subset \sum_{l}(A),$$

\noindent which means that $A$ is $l$-complete, impossible. Thus $ r
<Q^2$.

\noindent We define a new $A^{\flat}$ by taking the union of the
existing one with $B$. The bound on $|B|$ shows that  the new
$A^{\flat}$ is still of cardinality $O((pm)^{6/13}\log^2 p)$. We
keep using the notation $A'$ for $A\backslash A^{\flat}$, but the
reader should keep in mind that the new  $A'$ has no nice pair as we
have discarded $B$. This implies that  there are  intervals
$A_0,\dots,A_n$ of $\Z_p$ such that $|A_i| \le p/Q^2$ and $\min
\{\|x-y\| \Big| x \in A_i, y \in A_j \} \ge L$ for any $i \neq j$
and the union $\cup_{i=1}^n A_i$ contains $A'$.  It then follows
that

$$n+1 \le p/L.$$

\noindent But by pigeon-hole principle there is an interval, say
$A_0$, which contains at least $|A'|/(n+1)$ elements of $A'$. Recall
that the length of $A_0$ is less than $p/Q^2$ and

$$|A'|/(n+1) \ge |A'|L/p > (pm)^{6/13+12/13}/(pm)=(pm)^{5/13} > 3Q.$$

\end{proof}

\noindent We infer from Lemma \ref{lemma:main3:4} that, by an
appropriate translation, one can find a reasonably short interval
around 0 which contains many elements of $A$. (Notice that the
translation shifts $P$ to another interval of the same length). We
will work with this translated image of $A$.

\subsection{Distribution of the elements of $A$}

\noindent Let $I_0$ and $J_0$ be two disjoint subsequences of
$A'\cap [-p/(2Q^2),p/(2Q^2)] $ of cardinality $Q$ and $2Q$
respectively.

\noindent Let $A''=A'\backslash (I_0\cup J_0)$. We show that almost
all elements of $A''$( and thus almost all elements of $A$)
concentrate around a few points which are regularly distributed in
$\Z_p$.

\begin{lemma}\label{lemma:main3:1}

\noindent There is a subsequence $A'''\subset A''$ and an integer
$D$ such that

\begin{itemize}
\item $$|A'''|\le 2(pm)^{6/13},$$

\item $$D \le (pm)^{1/13},$$

\item for any $a\in A''\backslash A'''$ there is an integer
$0\le h \le D-1$ satisfying

$$|a-\frac{hp}{D}| \le \frac{p}{Q}.$$

\end{itemize}

\end{lemma}

\noindent We postpone the proof of Lemma \ref{lemma:main3:1} until
Proposition \ref{prop:main3:2}.

\noindent Let $a$ be any element of $A''$. Then by Dirichlet's
theorem, there is a pair of positive integers $i$ and $d$ satisfying
$1\le d\le Q$ and $\gcd(i,d)=1$ such that

$$|a - \frac{ip}{d}|\le \frac{p}{dQ}.$$

\noindent Next let

$$X_d=\{ a \in A'': |a - \frac{ip}{d}|\le
\frac{p}{dQ}, 1\le i \le d, 1\le d \le Q, \gcd(i,d)=1\}.$$

\noindent Call the index $d$ {\it rich} if $|X_d|\ge 2d.$ Let us
denote the {\it rich} indices by

$$d_1<d_2<\dots<d_s.$$

\noindent We will collect some facts about the {\it rich} indices.

\begin{proposition}\label{prop:main3:1}

$$d_j \le (pm)^{1/13}.$$

\end{proposition}

\begin{proof} (Proof of Proposition \ref{prop:main3:1}) Let
$X_{d_j}'=\{a_1,\dots, a_{d_j}\}$ be any subsequence of $d_j$
elements of $X_{d_j}$. By Lemma \ref{lemma:simple2}, for $0\le i \le
d_j-1$ there exists $A_{d_j}^i \subset X_{d_j}'$ such that

$$|\sum_{a \in A_{d_j}^i} a - \frac{ip}{d_j}| \le \frac{p}{Q}.$$

\noindent Choose a sequence $B_{d_j}^i \subset I_0$ such that
$|B_{d_j}^i|=d_j - |A_{d_j}^i|.$ By the definition of $I_0$ we have

$$\sum_{b \in B_{d_j}^i} |b| \le |B_{d_j}^i|p/(2Q^2) \le d_jp/(2Q^2) \le
p/2Q.$$

\noindent Thus

\begin{equation}\label{equation:main3:2}
|\sum_{a \in A_{d_j}^i} a + \sum_{b \in B_{d_j}^i} b -
\frac{ip}{d_j}| \le 2p/Q.
\end{equation}

\noindent By definition, $\sum_{a \in A_{d_j}^i} a + \sum_{b \in
B_{d_j}^i} b \subset \sum_{d_j}(X_{d_j} \cup I_0)$. Thus the
inequality (\ref{equation:main3:2}) implies that
$\sum_{d_j}(X_{d_j}' \cup I_0)$ forms a $K$-net of $\Z_p$ with $ K
\le p/d_j + 4p/Q$.

\noindent Now we claim that $K>L$. Seeking a contradiction, suppose
that $K\le L$. Then

\begin{equation}\label{equation:main3:3}
\sum_{d_j}(X_{d_j}' \cup I_0)+P=\Z_p.
\end{equation}

\noindent Because the cardinality of $A''\backslash X_{d_j}'$ is
larger than $l$,

$$|A''\backslash X_{d_j}'|= |A'|-|I_0|-|J_0|-|X_{d_j}'| \ge |A|-|A^\flat|-4Q \ge l,$$

\noindent we can choose $C \subset A''\backslash X_{d_j}'$ of
cardinality $|C|=l-d_j-l^\flat$. Next, by (\ref{equation:main3:3})
we have

$$\Z_p = \sum_{d_j}(X_{d_j}' \cup I_0) + P = \sum_{d_j}(X_{d_j}' \cup I_0) + P + \sum_{c\in C} c \subset
\sum_{l}(A).$$

\noindent Thus $A$ is $l$-complete, a contradiction.

\noindent Observe that beside the inequality $K > L$ we also have

$$L \gg p/Q \mbox{ and } L\ge \lambda(pm)^{12/13}/m \ge 2(pm)^{12/13}/m.$$

\noindent Thus

$$d_j \le 2p/L \le (pm)^{1/13}.$$

\end{proof}

\noindent Proposition \ref{prop:main3:1}, in particular, implies
that the number of {\it rich} indices is also small,

$$s\le (pm)^{1/13}.$$

\noindent In the following, we prove a stronger fact.

\begin{proposition}\label{prop:main3:2}

\noindent Let $D=lcm(d_1,\dots,d_s)$. Then we have

$$D \le (pm)^{1/13}.$$

\end{proposition}

\begin{proof} (Proof of Proposition \ref{prop:main3:2})
For each $1\le i \le s$ let $X_{d_i}'$ be a subsequence of
cardinality $d_i$ in $X_{d_i}$. We claim that $(\sum_{i=1}^s
d_i)^\ast (\bigcup_{i=1}^s X_{d_i}'\bigcup I_0)$ is a $K$-net in
$\Z_p$ with

$$K \le p/D + 4sp/Q.$$

\noindent To prove the claim, first let $r$ be any integer between 0
and $D-1$. By Lemma \ref{lemma:simple3} there exist $0\le a_i \le
d_i-1$ such that $\sum_{i=1}^s a_ip/d_i = rp/D.$

\noindent Next choose $A_{d_i}^r\subset X_{d_i}'$ such that
$|\sum_{a \in A_{d_i}^r} a -a_ip/d_i| \le p/Q$. Summing these
inequalities over $1\le i \le s$ we obtain

\begin{equation}\label{equation:main3:4}
|\sum_{a \in \bigcup_{i=1}^s A_{d_i}^r}a - rp/D| \le sp/Q.
\end{equation}

\noindent In addition, because

$$\sum_{i=1}^s d_i \le \lfloor s(pm)^{1/9} \rfloor \le \lfloor (pm)^{2/9} \rfloor =Q =|I_0|,$$

\noindent there are disjoint subsequences
$B_{d_1}^r,\dots,B_{d_s}^r$ of $I_0$ such that
$|B_{d_i}^r|=d_j-|A_{d_j}^r|$. And by the definotion of $I_0$ we
have

\begin{equation}\label{equation:main3:5}
\sum_{b \in \bigcup_{i=1}^s B_{d_i}^r} |b| \le (\sum_{i=1}^s
d_i)p/(2Q^2) \le Qp/(2Q^2) = p/2Q .
\end{equation}

\noindent Putting the estimates
(\ref{equation:main3:4}),(\ref{equation:main3:5}) together to obtain

\begin{equation}\label{equation:main3:6}
|\sum_{a \in \cup A_{d_i}^r}a + \sum_{b \in \cup B_{d_i}^r}b - rp/D|
\le sp/Q + p/2Q \le 2sp/Q.
\end{equation}

\noindent Notice that $\sum_{i=1}^s (|A_{d_i}^r| + |B_{d_i}^r|) =
\sum_{i=1}^s d_i$. Point (\ref{equation:main3:6}) concludes the
claim.

\noindent We now claim that $K > L$. Assume otherwise. Then

\begin{equation}\label{equation:main3:7}
\sum_{\sum_{i=1}^s d_i}(\bigcup_{i=1}^s X_{d_i}'\bigcup I_0) + P =
\Z_p.
\end{equation}

\noindent But
$$|A''\backslash \bigcup_{i=1}^s X_{d_i}'| =
|A'|-|I_0|-|J_0|-\sum_{j=1}^s d_j \ge |A|-|A^\flat|-4Q \ge l,$$

there exists a subsequence $C$ in $A''\backslash \bigcup_{i=1}^s
X_{d_i}'$ of cardinality $|C|=l-\sum_{j=1}^s d_j-l^\flat$.

\noindent Adding elements of $C$ to (\ref{equation:main3:7}) we
achieve

$$\Z_p = \sum_{d_j}(X_{d_j}' \cup I_0) + P = \sum_{d_j}(X_{d_j}' \cup I_0) + P + \sum_{c\in C} c .$$

\noindent The last sum of the equality above is a subset of
$\sum_l(A)$. Thus $A$ is $l$-complete, a contradiction.

\noindent In conclusion we have just proved that
$\sum_{d_1+\dots+d_s}(\bigcup_{i=1}^s X_{d_i}'\bigcup I_0)$ is a
$K$-net in $\Z_p$ with

$$L \le K \le p/D + 4sp/Q.$$

\noindent In particular,

$$L\le  p/D + 4sp/Q,$$

$$\lambda(pm)^{12/13}/m - 4p(pm)^{1/13}/(pm)^{3/13} \le p/D.$$

\noindent Hence (because $\lambda \ge 2$)

$$D\le (pm)^{1/13}.$$

\end{proof}

\noindent

For brevity set $t:=\sum_{i=1}^s d_i$, $H:=\bigcup_{i=1}^s
X_{d_i}'\cup I_0$ and

$$T:=\sum_t(H)=\sum_{d_1+\dots+d_s}(\bigcup_{i=1}^s X_{d_i}'\bigcup
I_0).$$

\noindent Recall that $T$ is a $K$-net with $K \le p/D + 4sp/Q$. We
remove $H$ from $A''$ and record the set $T$ for latter use. Let us
now prove Lemma \ref{lemma:main3:1} by putting everything together.

\begin{proof} (Proof of Lemma \ref{lemma:main3:1}) Call an element $a$ of $A''$
{\it single} if $a\not\in \bigcup_{j=1}^s X_{d_j}$. By Dirichlet's
theorem, any single point is an element of some $X_d$ where $d$ is
not rich. But $|X_d|<2d$ if $d$ is not rich. Thus by double
counting, the number of single points, denoted by $A'''$, is bounded
by

$$|A'''| \le \sum_{d\le Q}(2d-1) < 2Q^2=2(pm)^{6/13}.$$

\noindent Let $a$ be any element of $A''\backslash A'''$, then $a\in
X_{d_j}$ for some rich $d_j$. Put $h=iD/d_j$. Then by definition

$$|a-\frac{hp}{D}|=|a-\frac{ip}{d_j}|\le \frac{p}{d_j} \le  \frac{p}{Q}.$$

\noindent Furthermore, by Proposition \ref{prop:main3:1},

$$D\le (pm)^{1/13}.$$

\end{proof}

\noindent Add $A'''$ to $A^\flat$, the cardinality of $A^\flat$ is
still $O((pm)^{6/13}\log^2 p)$. For $1\le h \le D$ we let

$$J_h=\{a|a\in A'',  \frac{hp}{D}-\frac{p}{Q}\le a\le \frac{hp}{D} + \frac{p}{Q}\}$$

and

$$R_h= \{a-\frac{hp}{D}|a\in J_h\}.$$

\noindent By throwing away a small number ($\le sD\le (pm)^{2/13}$)
of elements to $A^{\flat}$, we can assume that the cardinalities of
$R_h$, $1\le h \le s$, are divisible by $D$. Note that the sum of
any $D$ elements of $R_h$ is an integer. We denote by $R$ the
sequence of all reduced elements,

$$R=\bigcup_{h=1}^s R_h. $$

\noindent Hence for any $r\in R$ we have $|r|\le p/Q$.

\noindent Let us summarize what we have obtained up to this step. Up
to a proper dilation (with $b'$) and translation (with $c'$), there
is a partition of $A$, $A=A^\flat \cup J_0 \cup H \cup A''$ such
that

\begin{itemize}

\item $|A^{\flat}|=O((pm)^{6/13}\log^2 p)$ and $\sum_{l^{\flat}}(A^{\flat})$ contains an
interval $P = [a,a+L]$ of length $L=\lambda (pm)^{12/13}/m$ with
some $l^\flat \le (pm)^{6/13}\log^2 p$ .

\vskip .1in

\item $|H| \le 2(pm)^{3/13}$ and $\sum_t(H)$ contains a $p/D + 4sp/Q$-net
(named $T$).

\vskip .1in

\item $|J_0|=2Q $ and $J_0\subset [-p/(2Q^2), p/(2Q^2)].$

\end{itemize}

\subsection{Completing the proof of Theorem \ref{theorem:main3}} Set
$$l_0:=l-l^\flat -t.$$
Since the elements of $R$ are small, the set $\sum_{l_0}R$ (which is
a subset of $\Q$ rational numbers) is dense in the interval in which
it is contained. We show that $\sum_{l_0}(R) \cap \Z$ is also dense
in this interval. Suppose for the moment that this interval is
longer than $p/D + 4sp/Q$. Then $(\sum_{l_0}(R) \cap \Z) + P$
contains another interval of length $p/D + 4sp/Q$ (in $\Z$, as $P$
is viewed as an interval of $\Z$). We then infer that
$(\sum_{l_0}(A'') \cap \Z) + P$ contains an interval of that same
length in $\Z_p$. So

$$\Z_p = \sum_{l_0}(A'') + P + T \subset \sum_{l}(A).$$

\noindent Which is impossible. We conclude that $\sum_{l_0}(R)$ must
be supported by a short interval of $\Q$. In the following we
explain the argument in detail.

\noindent Set

$$l_2:=l_0 - D^2 \mbox{ and } l_1:=l_2-Q=l_2-\lfloor (pm)^{1/3}\rfloor.$$

\noindent Then

$$l_2 > l_1 \ge l- 2(pm)^{6/13}\log^2 p.$$

\noindent Viewing $R$ as a subsequence of $\Q$ in $[-p/Q,p/Q]$, our
goal is to establish the following.

\begin{lemma}\label{lemma:main3:2}
Let $m_1=\min_{l_1\le l'\le l_2} (\min(\sum_{l'}(R))$ and
$m_2=\max_{l_1\le l'\le l_2} (\max \sum_{l'}(R)) $. Then we have

$$m_2-m_1 < p/D.$$

\end{lemma}

\begin{proof} (Sketch of proof of Lemma \ref{lemma:main3:2})
Add several (at most $D^2$) elements of $R$ to the representations
of $m_1$ and $m_2$ respectively to make the number of summands from
each class $R_h$ divisible by $D$. We obtain $m_1',m_2'$ with the
following properties.

\begin{itemize}

\item $m_i' \in \sum_{l_i'}(R)$, where $l_1 \le l_i'\le l_2+D^2.$

\vskip .1in

\item $|m_i'-m_i| \le D^2p/Q.$ (Because to create $m_i'$ we added at most $D^2$
elements from $R$, whose element is bounded by $p/Q$.)

\vskip .1in

\item $m_1',m_2'\in \Z.$ (As the sum of any $D$ elements of
$R_h$ is an integer.)

\end{itemize}

\noindent By the properties above, we are done with the Lemma if
$m_2'-m_1' \le p/D-2D^2p/Q$. Seeking for contradiction, suppose that

\begin{equation}\label{equation:main3:9}
m_2'-m_1' > p/D-2D^2p/Q.
\end{equation}

\noindent Let $U_1,U_2\subset R$ be sequences of cardinality
$l_1',l_2'$ respectively such that

$$\sum_{u \in U_i} u= m_i'.$$

\noindent The reader should find it straightforward to construct
sequences $V_1,V_2, \dots, V_n$ in $R$ such that all the following
properties hold.

\begin{itemize}

\item $$V_1=U_1, V_n=U_2.$$ \\

\item $$\min \{l_1',l_2'\} \le |V_i| \le \max \{l_1',l_2'\} \mbox{ for $1\le i \le n$}.$$\\

\item \begin{equation}\label{equation:main3:11} |V_{i+1}\backslash
V_i|\le D.
\end{equation}\\

\item For any $1\le h \le s$ the cardinality of $V_i\cap R_h$ is divisible by
$D$, i.e.,

\begin{equation}\label{equation:main3:10}
D||V_i\cap R_h| \mbox{ for } 1\le h \le s.
\end{equation}

\end{itemize}

\noindent Notice that condition (\ref{equation:main3:10}) guarantees
that $\sum_{v\in V_i}v$ is an integer, and (\ref{equation:main3:11})
implies that

$$|\sum_{v \in V_{i+1}}v - \sum_{v\in V_i} v| \le Dp/Q \mbox{ for } 1\le i \le n.$$

\noindent Thus the set $\{\sum_{v \in V_i} v| i=1,\dots,n \}$ is a
$pD/Q$-net (of $\Z$) in the interval $[m_1',m_2']$. Recall that

$$|J_0|=2Q > Q+D^2 = l_0 - l_1 \ge l_0 - |V_i|,$$

\noindent i.e. for each $1\le i\le n$ one can choose a sequence
$W_i$ of $J_0$ of cardinality $l_0-|V_i|$ ($W_i$'s are not
necessarily disjoint). Denote $V_i\cup W_i$ by $X_i$. Then we have
$|X_i|=l_0$ and

\begin{equation}\label{equation:main3:12}
|\sum_{x\in X_i} x - \sum_{v\in V_i} v|\le (l_0-|V_i|)p/Q^2 \le
(l_0-l_1)p/Q^2 \le p/Q.
\end{equation}

\noindent Because $\{\sum_{v \in V_i} v| i=1,\dots,n \}$ is a
$Dp/Q$-net in $[m_1',m_2']$, we have

$$[m_1',m_2'] \subset \{\sum_{v \in V_i} v| i=1,\dots,n \}+ [0,Dp/Q]
(\mod p);$$

\noindent and it follows from (\ref{equation:main3:12}) that

\begin{equation}\label{equation:main3:13}
[m_1'+p/Q,m_2'-p/Q] \subset \{\sum_{x\in X_i} x |i=1,\dots,n
\}+[0,2Dp/Q].
\end{equation}

\noindent We proceed by claiming the following.

\begin{claim}\label{fact:main3:1}
Suppose that (\ref{equation:main3:9}) holds. Then the set
$$\{\sum_{x\in X_i} x + T|i=1,\dots,n \}$$
is a $8D^2p/Q$-net of $\Z_p$.
\end{claim}

\begin{proof} (Proof of Claim \ref{fact:main3:1}) Obtain from (\ref{equation:main3:13})
that

$$[m_1'+p/Q,m_2'+ 7D^2p/Q ] \subset \{\sum_{x\in X_i} x |i=1,\dots,n
\}+[0,8D^2p/Q].$$

Consequently,
\begin{equation}\label{equation:main3:14}
[m_1'+p/Q,m_2'+ 7D^2p/Q ]+T \subset \{\sum_{x\in X_i} x |i=1,\dots,n
\}+[0,8D^2p/Q] + T.
\end{equation}

\noindent Notice that because $T$ is a $p/D + 4sp/Q$-net of $\Z_p$,
and by (\ref{equation:main3:9}) that

$$m_2'+ 7D^2p/Q- m_1'- p/Q \ge p/D+ 4D^2p/Q > p/D + 4sp/Q,$$

we have

$$\Z_p=[m_1'+p/Q,m_2'+ 7D^2p/Q ]+T.$$

\noindent Together with (\ref{equation:main3:14}) this gives

$$(\{\sum_{x\in X_i} x |i=1,\dots,n \}+ T) + [0,8D^2p/Q] = \Z_p.$$

\end{proof}

\noindent To finish the proof of Lemma \ref{lemma:main3:2} one
observes that

$$L \ge \lambda p/(pm)^{1/13} \ge 8p/(pm)^{1/13} \ge 8D^2p/Q.$$

\noindent Thus Claim \ref{fact:main3:1} would give

$$\{\sum_{x\in X_i} x + T|i=1,\dots,n \} + P = \Z_p.$$

\noindent However, $\{\sum_{x\in X_i} x + T|i=1,\dots,n \} + P
\subset \sum_{l}(A)$. Hence $A$ is $l-$complete, a contradiction. As
a consequence, (\ref{equation:main3:9}) can not hold.

\end{proof}

\noindent Now we close the proof of Theorem \ref{theorem:main3}.
Dilate the whole set $A$ with $D$. By viewing $D \cdot A''$ as a
sequence of $\Z$ in $[-Dp/Q,Dp/Q]$, one sees that

$$\max_{l_1 \le l' \le l_2}\max(\sum_{l'}(D \cdot A'')) - \min_{l_1\le l'\le l_2}\min(\sum_{l'}(D \cdot A'')) = Dm_2 -Dm_1 < p.$$

\noindent Thus if $\Phi$ denotes the linear map $b' \cdot X+c'$ then
the statement of Theorem \ref{theorem:main3} holds for $A^\flat
(\mbox{ of the statement }):=\Phi^{-1}(A^{\flat} \cup J_0 \cup H)$
and $b:=Db',c:=Dc'$.

\section{Sketch of proof of Theorem
\ref{theorem:main2}}\label{section:main2}

\noindent Theorem \ref{theorem:main2} can be verified by following
the proof of Theorem \ref{theorem:main3} above. In fact, the
situation here is somewhat simpler. Since the subset sums in Theorem
\ref{theorem:main2} do not need to have a fixed number of summands,
we do not have to consider $I_0$ and $J_0$.

\noindent Keep the same notation as in the proof of Theorem
\ref{theorem:main3}. As an analogue of  Lemma \ref{lemma:main3:2},
we can establish the following lemma.

\begin{lemma}\label{lemma:main2:1}
Let $m_1=\min(\sum(R))$ and $m_2=\max(\sum(R)) $. Then we have

$$m_2-m_1 < p/D.$$

\end{lemma}

\noindent Then by dilating the whole set $A$ with $D$, one obtains
Theorem \ref{theorem:main2}.

\section{Proof of Theorem
\ref{theorem:main1}}\label{section:main1}

\noindent By Theorem \ref{theorem:main2} there exists a non-zero
residue $b$ and a small set $A^{\flat}\subset A$ of cardinality at
most $c_2f(p,m)$ such that

\begin{equation}\label{equation:main1:1}
\sum_{a\in b\cdot (A\backslash A^{\flat})} \|a\| <p.
\end{equation}

\noindent Consider the sequence of positive and negative elements of
$b\cdot (A\backslash A^{\flat})$,

$$A^+:=b \cdot (A\backslash A^{\flat})\cap [1,(p-1)/2] \mbox{ and } A^-:=
b\cdot (A\backslash A^{\flat})\cap [-(p-1)/2,-1].$$

\noindent We shall prove the following.

\begin{lemma}\label{lemma:main1:1}
There exists an absolute constant $\beta$ such that either  $|A^+| \le \beta f(p,m)$ or  $|A^-| \le \beta f(p,m).$
\end{lemma}

\noindent Assume for the moment, and without loss of generality,
that $|A^-| \le \beta f(p,m)$. Then one may verify that Theorem
\ref{theorem:main1} holds for $A^\flat \mbox{ (of Theorem
\ref{theorem:main1}) }:= A^\flat \cup b^{-1} \cdot A^-$ and
$c_1:=c_2+\beta$. Thus it remains to prove Lemma
\ref{lemma:main1:1}.

\begin{proof} (Proof of Lemma \ref{lemma:main1:1}) Assume otherwise that

\begin{equation}\label{equation:main1:2}
|A^+|,|A^-| \ge \beta f(p,m) \mbox { for large positive
constant } \beta.
\end{equation}

Note that from (\ref{equation:main1:1}) we have

\begin{equation}\label{equation:main1:3,4}
\sum_{a\in A^+}a < p, \mbox{ and } \sum_{a\in A^-}|a| < p.
\end{equation}

Set $q:=\lfloor p/f(p,m)\rfloor$. Let $B^+:=A^+\cap [1,q]$ and
$B^-:=A^-\cap [-1,-q]$ respectively.

\noindent We infer from (\ref{equation:main1:3,4}) that

$$|B^+|\ge (\beta-1)f(p,m) \mbox{ and } |B^-| \ge (\beta-1)f(p,m).$$

\noindent Viewing $B^+$ and $B^-$ as sequence of integers in $[-q,q]$,
we then reach a contradiction with the zero-sum-freeness property of $A$ by showing that there exist some elements of $B^+$ and $B^-$ whose sum is 0.

\noindent Consider the following two cases.

{\bf Case 1:} $m\ge p^{4/9}$.

\noindent By pigeon-hole principle there are two elements $a^+ \in
B^+,a^-\in B^-$ whose multiplicities (denoted by $m_{a^+},m_{a^-}$
respectively) are large.

$$m_{a^+} \ge |B^+|/q \ge (\beta-1)f(p,m)/(p/f(p,m)) > (pm)^{12/13}/p > p/(pm)^{6/13} \ge q,$$

\noindent and similarly

$$m_{a^-} > q.$$

\noindent Note that $0\le |a^-|,a^+ \le q$. Thus $|a^-| < m_{a^-}$
and $a^+ < m_{a^+}$, which yield

$$0=|a^-|a^+ + a^+ a^- \in S_{B^+} + S_{B^-} \subset \sum(A), \mbox{contradiction}.$$

{\bf Case 2:} $1\le m <p^{4/9}$.

\noindent Without loss of generality assume that

\begin{equation}\label{equation:main1:5}
|\sum_{a\in B^-} a | \ge \sum_{a\in B^+} a.
\end{equation}

\noindent Fix any subset $X$ of $B^+$ of cardinality
$|X|=f(p,m)/\max(\log{p},m)$.

\noindent First, one sees that

$$(f(p,m)/\log p)^2 \gg  p/f(pm)=q,$$

and

$$(f(p,m)/m)^2 \gg p/f(p,m)=q.$$

\noindent Thus, Theorem \ref{lemma:main1} applied to $X$ (with
$l=\lfloor |X|/2 \rfloor$ and $d=1$) yields an arithmetic
progression $P=\{a,a+d,\dots, a+Ld\}$ of length $L \ge c(1)|X|^2/2$.

\noindent Note that $P\subset S_{X}\subset [1,|X|q ]$, thus the
difference $d$ of $P$ is bounded, i.e.,

\begin{equation}\label{equation:main1:6}
d \le|X|q/L \le 2q/(c(1)|X|) \ll (pm)^{1/13}/\log{p}.
\end{equation}

\noindent Next, view $(B^+\backslash X)\cup B^-$ as a sequence of
residues modulo $d$. We throw away residues of multiplicity less
than $d$. Let $W$ be the sequence of thrown elements. So obviously,

$$|W| \le d^2 \le (pm)^{2/9}/\log^2{p}.$$

\noindent We consider two subcases.

{\bf Subcase 2.1:} There exists a nontrivial divisor $d_1$ of $d$
which divides all the remaining residues.

\noindent Set

$$B^+_1:=\{\frac{b}{d_1}|b\in B^+\backslash (X\cup
W)\} \mbox{ and } B^-_1:=\{\frac{b}{d_1}| b\in B^-\backslash W \}.$$

\noindent Observe that

$$|B^+_1|,|B^-_1| \ge (\beta-1) f(pm)-2f(p,m)/\log p.$$

\noindent Also, $B^+_1,B^-_1 \subset [-q_1,q_1]:=[-\lfloor q/d_1
\rfloor , \lfloor q/d_1 \rfloor ]$.

\noindent Viewing $B^+_1$ and $B^-_1$ as $B^+$ and $B^-$, we
reconsider {\bf Case 1} and {\bf Case 2}. Thus either a
contradiction is obtained or we get $B^+_2$ and $B^-_2$ whose
elements are divisible by some integer $d_2 \ge 2$. Repeat the
process until we get a contradiction thanks to {\bf Case 1} or {\bf
Subcase 2.2} as follows. (Notice that the process stops after at
most $\log p$ steps because $q_i$ decreases by a factor of at least
2 with each step, while $|B^+_i|,|B^-_i|\ge (\beta -2)f(p,m)$
always.)

\vskip 2mm

{\bf Subcase 2.2:} There does not exist such divisor of $d$. Thus
the residues are mutually co-prime with $d$.

\noindent By Lemma \ref{lemma:simple4} there exist $x_1, \dots, x_u
\in X \backslash X, y_1,\dots,y_v \in B^-$ with $u+v\le d$ and

\begin{equation}\label{equation:main1:7}
a=-\sum_{i=1}^u x_i - \sum_{j=1}^v y_j  (\mod d).
\end{equation}

\noindent Note that

\begin{equation}\label{equation:main1:8}
\sum_{i=1}^u x_i + |\sum_{j=1}^v y_j|\le dq \ll dL.
\end{equation}

\noindent We consider the following two possibilities.

{\bf Subcase 2.2.1:} $|\sum_{j=1}^v y_i| - \sum_{i=1}^u x_i \ge a.$

\noindent Then by (\ref{equation:main1:7}) and
(\ref{equation:main1:8}) we get

$$|\sum_{j=1}^v y_i| - \sum_{i=1}^u x_i \in P.$$

\noindent Thus

$$|\sum_{j=1}^v y_i| \in \sum_{i=1}^u x_i + \sum(X),$$

and so

$$0 \in \sum_{j=1}^v y_i + \sum_{i=1}^u x_i + \sum(X) \subset \sum(X
\cup B^-) \subset \sum(A), \mbox{ contracdition }. $$

\vskip 5mm

{\bf Subcase 2.2.2:} $|\sum_{j=1}^v y_i| - \sum_{i=1}^u x_i < a.$

\noindent Then let $Y_0=:\{y_1,\dots,y_v\}$. By Lemma
\ref{lemma:simple2} one can find $Y_0' \subset B^- \backslash Y_0$
such that $|Y_0'|\le d$ and $d | \sum_{y \in Y_0'} y.$

\noindent Set $Y_1:=Y_0 \cup Y_0'$. If $|\sum_{y \in Y_1} y | -
\sum_{i=1}^u x_i$ is still less than $a$ then we again use Lemma
\ref{lemma:simple2} to find $Y_1' \subset B^-\backslash Y_1$ such
that $Y_1'$ has the same property as $Y_0'$. We next increase $Y_1$
by $Y_2:=Y_1\cup Y_1'$. Repeat the process until we get $Y_N \subset
B^-$ such that

$$|\sum_{y \in Y_{N-1}} y | - \sum_{i=1}^u x_i < a \mbox{ and } |\sum_{y \in Y_N} y | - \sum_{i=1}^u x_i \ge a.$$

\noindent Notice that by (\ref{equation:main1:3,4}) we have

$$\sum_{i=1}^u x_i + a + Ld \le \sum_{a\in B^+} a \le |\sum_{y\in
B^-} y|.$$

\noindent In addition, since $q \ll L$,

\begin{equation}\label{equation:main1:9}
\sum_{i=1}^u x_i + a \le |\sum_{y\in Y} y|
\end{equation}

\noindent for any $Y\subset B^-$ with cardinality $|Y|\ge |B^-|-d$.
Lemma \ref{lemma:simple2} and (\ref{equation:main1:9}) thus ensure the
existence of $N$ above.

\noindent In sum,
$$0 \le |\sum_{y\in Y_N} y |-\sum_{i=1}^u x_i-a \le Ld$$

\noindent and $d$ is divisible by $|\sum_{y\in Y_N} y |-\sum_{i=1}^u
x_i-a.$

\noindent It then follows that

$$|\sum_{y\in Y_N} y |-\sum_{i=1}^u x_i \in P.$$

$$0 \in \sum_{y \in Y_N} y + \sum_{i=1}^u x_i + \sum(X') \subset \sum(X
\cup Y) \subset \sum(A), \mbox{ contradiction }.$$

\end{proof}


\begin{thebibliography}{99}

\bibitem{A} G. E. Andrews, {\it The theory of partitions.} Cambridge
university press, 1998.



\bibitem{CS} {F. Chen} and {S. Savchev}, {\it Long n-zero-free
sequences in finite cyclic groups}. Discrete Mathematics, 308,
(2008), 1-8.


\bibitem{D1} { J. M. Deshouillers}, {\it Quand seule la sous-somme
vide est nulle modulo $p$}, the prodeeding of the Journees
Arithmetiques 2005.

\bibitem{D2} { J. M. Deshouillers}, {\it Lower bound concerning subset sum
wich do not cover all the residues modulo $p$}, Hardy- Ramanujan
Journal, Vol. 28(2005) 30-34.

\bibitem{DF} { J. M. Deshouillers} and {Gregory A. Freiman}, {\it When
subset-sums do not cover all the residues modulo $p$}, Journal of
Number Theory 104(2004) 255-262.


\bibitem{EGZ} { P. Erd\H{o}s,
A. Ginzburg} and {A. Ziv}, {\it Theorem in the additive number
theory}. Bull. Res. Council Israel 10F (1961), 41-43.


\bibitem{GG} {W. D. Gao and A. Geroldinger}, {\it Zero-sum problems in finite abelian
groups : a survey}, Expo. Math. 24 (2006), 337-369.


\bibitem{GLP1}{W. D. Gao, Y. Li, J. Peng} and {F. Sun}, {\it Subsums of a zero-sum free subset of an abelian group},
E. J. Comb. 15 (2008), Research Paper 116.

\bibitem{GLP2} {W. D. Gao, Y. Li, J. Peng} and {F. Sun}, {\it On subsequence sums of a
zero-sum free sequence II}, E. J. Comb. 15 (2008), Research paper
117.

\bibitem{GPT} { W. D. Gao, A. Panigrahi} and {R. Thangadurai},
{\it On the structure of $p$-zero-sum-free sequences and its
application to a variant of Erd\H{os}-Ginzburg-Ziv theorem}. Proc.
Indian Acad. Sci. Vol. 115, No. 1 (2005), 67-77.


\bibitem{Ge} {A. Geroldinger}, {\it Additive group theory and non-unique
factorizations, Combinatorial Number Theory and Additive Group
Theory}, Advanced Courses in Mathematics CRM Barcelona,
Birkh\"auser, 2008.


\bibitem{NgSzemVu} { H. H. Nguyen, E. Szemer\'edi} and  {V. H. Vu},
{\it Subset sums modulo a prime}, Acta Arithmetica, 131.4 (2008),
303-316.


\bibitem{Olson2}  { J. E. Olson}, {\it Sums of sets of group elements}.
Acta Arithmetica, 28 (1975), 147-156.

\bibitem{Sun}  {J. W. Sun}, {\it List of publications on restricted sumsets}.
2005.

\bibitem{SzemVu1} { E. Szemer\'edi} and { V. H. Vu },
{\it Long arithmetic progression in sumsets and the number of x-free
sets}. Proceeding of London Math Society, 90(2005)
273-296.

\bibitem{SzemVu2} { E. Szemer\'edi} and { V. H. Vu }, {\it
Long arithmetic progressions in sumsets: Thresholds and Bounds}.
Journal of the A.M.S, 19 (2006), no 1, 119-169.

\end{thebibliography}
\end{document}